\documentclass{article}

\usepackage{amsmath}

\begin{document}

\title{Some Applications of the Method of Normal\\ Fundamental Functions
to Oscillation Problems\footnote{Accepted (21-March-2006) to the
Proceedings of MTNS 2006 (the 17th International Symposium on
Mathematical Theory of Networks and Systems), to be held on July
24-28, 2006, in Kyoto, Japan. Research Report CM06/I-17,
University of Aveiro, May 2006.}}

\author{Olena V. Mul\thanks{Supported by FCT (the Portuguese Foundation for Science
and Technology), fellowship SFRH/BPD/14946/2004.}\\
        \texttt{olena@mat.ua.pt}
        \and
        Delfim F.~M.~Torres\\
        \texttt{delfim@mat.ua.pt}}

\date{Department of Mathematics\\
      University of Aveiro\\
      3810-193 Aveiro, Portugal}

\maketitle


\begin{abstract}
We report on the possibilities of using the method of normal
fundamental systems for solving some problems of oscillation
theory. Large elastic dynamical systems with continuous and discrete
parameters are considered, which have many different engineering
applications. Intensive oscillations in such systems are possible,
but not desirable. Therefore, it is very important to obtain
conditions for which oscillations take or not-take place.
Mathematically, one needs to search for the solutions of partial
differential equations satisfying both boundary and conjugation
conditions. In this paper we overview the methodology of normal
fundamental systems for the study of such oscillation problems,
which provide an efficient and reliable computational method. The
obtained results permit to analyze the influence of different system
parameters on oscillations as well as to compute the optimal feedback
parameters for the active vibration control of the systems.
\end{abstract}


\smallskip


\smallskip


\smallskip


\section{Introduction}

Boundary problems with continuous and discrete parameters constitute
a wide class of oscillation problems with a great relevance in
mathematics and applications. They appear, for instance, in the
study of mechanical oscillations of ship hulls, ship masts,
antennas, supporting surfaces of aircrafts, turbine fins and shafts,
structural components of automatic control apparatus, etc. The
structural study of the dynamics for such elastic systems imply
solving some partial differential equations with non-stationary
continuous and discrete coefficients. Moreover, the solutions of the
partial differential equations have to satisfy not only boundary
conditions, but also some conjugation constraints, which
considerably complicate the problem. We claim that the method of the
normal fundamental system of solutions is a very efficient and
reliable way for numerically solving such boundary problems with
continuous and discrete parameters.


\section{The Method}

Many important problems of oscillation theory may be reduced to the
study of systems of partial differential equations with variable
coefficients, which describe oscillation processes in very different
systems and which may be written in the form:
\begin{equation}
\label{eq:1}
\begin{gathered}
\frac{{\partial z_k }}{{\partial y}}
= \sum\limits_{j=1}^N A_{kj} (y)z_j
+  \sum\limits_{j = 1}^N B_{kj}(y)\frac{{\partial ^2 z_j }}{{\partial t^2 }}
+  \sum\limits_{j = 1}^N C_{kj} (y)\frac{{\partial z_j }}{{\partial t}} \, , \\
y_0  \le y \le y_n \, , \quad k = \overline{1,N} \, ,
\end{gathered}
\end{equation}
where $z_j(y,t)$ is the deviation function of the oscillation system
from an equilibrium position, and the coefficients $A_{kj}(y)$,
$B_{kj}(y)$ and $C_{kj}(y)$ are real-valued piecewise continuous
functions bounded on the interval $[y_0 ,y_n ]$, \textrm{i.e.},
$\left| {A_{kj} (y)} \right| \le M$, $\left| {B_{kj} (y)} \right|
\le M$, $\left| {C_{kj} (y)} \right| \le M$ for some constant $M$
and $k,j = \overline{1,N}$.

Besides, we have the next linear homogeneous boundary conditions:
\begin{equation}
\label{eq:2}
\begin{gathered}
z_r(0) = 0 \, , \quad r = \overline{1,m} \, , \\
z_s (l) = 0 \, , \quad s = \overline{m + 1,N} \, ,
\end{gathered}
\end{equation}
where from here always $N = 2m$.

We assume that the conjugation conditions take place at the given
arbitrary discontinuity points $y = y_i$ of functions $A_{kj}(y)$,
$B_{kj}(y)$ and $C_{kj}(y)$, and they may be presented as follows:
\begin{multline}
\label{eq:3} \sum\limits_{k = 1}^N a_{k,p}^{(i + 1)} (y_i )z_k^{(i +
1)} (y_i ) +  \sum\limits_{k = 1}^N {b_{k,p}^{(i + 1)} (y_i)
\frac{{\partial ^2 z_k^{(i + 1)} (y_i )}}{{\partial t^2 }}} \\
+ \sum\limits_{k = 1}^N d_{k,p}^{(i + 1)}(y_i )\frac{{\partial
z_k^{(i + 1)} (y_i )}}{{\partial t}}
=  \sum\limits_{k = 1}^N {a_{k,p}^{(i)} (y_i )z_k^{(i)} (y_i )} \\
+ \sum\limits_{k = 1}^N {b_{k,p}^{(i)} (y_i) \frac{{\partial ^2
z_k^{(i)}(y_i )}}{{\partial t^2 }}} + \sum\limits_{k = 1}^N
{d_{k,p}^{(i)}(y_i )\frac{{\partial z_k^{(i)} (y_i )}}{{\partial
t}}} \, ,
\end{multline}
where $p = \overline{1,N}$, $i = \overline{1,n - 1}$ and the
superscript $i$ means the value of the corresponding function on the
$i$-th interval. In many practical applications conditions
\eqref{eq:3} are not taken into account because of significant
complication of investigations and this causes qualitative and
quantitative errors in the results \cite{c2}.

The continuous-discrete boundary problem \eqref{eq:1}-\eqref{eq:3}
serves as a model for a large number of oscillation dynamical
systems \cite{c2}.

The standard way for integration of system \eqref{eq:1}-\eqref{eq:3}
consists in applying the Fourier method of separation of variables
\cite{Kach-teor}. The solution then is searched in the form of a
product of two functions, where the first one depends only on the
coordinate $y$, and the second one depends only on time $t$:
\begin{equation*}
\label{eq:4}
z_k (y,t) = \varphi_k (y)e^{\lambda t} \, .
\end{equation*}
The problem is thus reduced to the one of finding the values for the
parameter $\lambda$, at which there exist non-trivial solutions
$\varphi_k (y)$ of the system \eqref{eq:1}-\eqref{eq:3}. Such values
$\lambda$ are called eigenvalues of the boundary problem, and the
corresponding solutions $\varphi_k (y)$ are called eigenfunctions.
The eigenvalues for the oscillation system, described by equations
\eqref{eq:1}-\eqref{eq:3}, are complex quantities $\lambda  = q +
ip$. The oscillations have increasing or decreasing amplitude,
depending on the sign of the real part $q$ of the eigenvalue.

In many practical oscillation theory problems, one is interested to
know eigenvalues at which oscillations with a constant amplitude
occur. That is only possible at an imaginary eigenvalue $\lambda =
ip$, where $p$ is the eigenfrequency of the oscillations of the
considered system. We assume that the system makes harmonic
oscillations with some constant amplitude, searching for harmonic
functions $z_k(y,t)$ of time. In this case we have $\varphi_k (y) =
\varphi _k^{(1)} (y) + i\varphi _k^{(2)}(y)$, where $\varphi_k^{(1)}
(y)$ and $\varphi _k^{(2)}(y)$ are real-valued functions. Hence, the
next form of solution may be used:
\begin{equation}
\label{eq:5} z_k(y,t) = [\varphi _k^{(1)} (y) + i \varphi _k^{(2)}
(y)]e^{ipt} \, .
\end{equation}

In particular cases, one of the functions $\varphi_k^{(1)} (y)$ or
$\varphi _k^{(2)}(y)$ may be identically equal to zero. Thus, if
$\varphi _k^{(2)}(y)\equiv 0$, then $\varphi _k(y)$ is a real-valued
function. This happens, for instance, in oscillatory motions of
mechanical systems with no forces proportional to the velocity, i.e.
in the absence of the terms with derivatives $\displaystyle \frac{\partial z_{j}}{\partial t}$.

After substitution of \eqref{eq:5} into equations \eqref{eq:1} and
separation of imaginary and real parts, we obtain respectively:

\begin{multline}
\label{eq:5a} \frac{{d\varphi_k^{(1)}}}{{dy}} = \sum\limits_{j =
1}^N {A_{kj} (y) \varphi_j^{(1)}}(y)- p^{2}\sum\limits_{j = 1}^N
{B_{kj} (y) \varphi_j^{(1)}}(y) - p \sum\limits_{j = 1}^N {C_{kj}
(y) \varphi_j^{(2)}}(y) \, , \\
\frac{{d\varphi_k^{(2)}}}{{dy}} = \sum\limits_{j = 1}^N {A_{kj} (y)
\varphi_j^{(2)}}(y)- p^{2}\sum\limits_{j = 1}^N {B_{kj} (y)
\varphi_j^{(2)}}(y) + p \sum\limits_{j = 1}^N {C_{kj} (y)
\varphi_j^{(1)}}(y) \, ,
\end{multline}
where $k = \overline{1,N}$.

The boundary conditions \eqref{eq:2} after similar transformations
may be written as

\begin{equation}
\label{eq:5b}
\begin{gathered}
\varphi_r^{(1)}(y_{0}) = \varphi_r^{(2)}(y_{0}) = 0 \, , \quad r = \overline{1,m} \, , \\
\varphi_s^{(1)}(y_{n}) = \varphi_s^{(2)}(y_{n}) = 0 \, , \quad s =
\overline{m + 1,N} \, .
\end{gathered}
\end{equation}

Analogously, the conjugation conditions \eqref{eq:3} will be
transformed to the next ones:

\begin{multline}
\label{eq:5c} \sum\limits_{k = 1}^N {a_{k,p}^{(i+1)} (y_{i})
\varphi_{k,i+1}^{(1)}}(y_{i})- p^{2}\sum\limits_{k = 1}^N
{b_{k,p}^{(i+1)} (y_{i}) \varphi_{k,i+1}^{(1)}}(y_{i}) \\
- p\sum\limits_{k = 1}^N {d_{k,p}^{(i+1)} (y_{i})
\varphi_{k,i+1}^{(2)}}(y_{i})= \sum\limits_{k = 1}^N
{a_{k,p}^{(i)} (y_{i}) \varphi_{k,i}^{(1)}}(y_{i}) \\
-p^{2}\sum\limits_{k = 1}^N {b_{k,p}^{(i)} (y_{i})
\varphi_{k,i}^{(1)}}(y_{i}) - p\sum\limits_{k = 1}^N {d_{k,p}^{(i)}
(y_{i}) \varphi_{k,i}^{(2)}}(y_{i})\, , \\
\sum\limits_{k = 1}^N {a_{k,p}^{(i+1)} (y_{i})
\varphi_{k,i+1}^{(2)}}(y_{i})- p^{2}\sum\limits_{k = 1}^N
{b_{k,p}^{(i+1)} (y_{i}) \varphi_{k,i+1}^{(2)}}(y_{i}) \\
+ p\sum\limits_{k = 1}^N {d_{k,p}^{(i+1)} (y_{i})
\varphi_{k,i+1}^{(1)}}(y_{i})= \sum\limits_{k = 1}^N
{a_{k,p}^{(i)} (y_{i}) \varphi_{k,i}^{(2)}}(y_{i}) \\
-p^{2}\sum\limits_{k = 1}^N {b_{k,p}^{(i)} (y_{i})
\varphi_{k,i}^{(2)}}(y_{i}) + p\sum\limits_{k = 1}^N {d_{k,p}^{(i)}
(y_{i}) \varphi_{k,i}^{(1)}}(y_{i})\, ,
\end{multline}
where $p = \overline{1,N}$.

Since the initial problem \eqref{eq:1}--\eqref{eq:3} is reduced to
the form \eqref{eq:5a}--\eqref{eq:5c}, we can consider now the
boundary problem for the system of ordinary differential equations,
or, more precisely, the problem with the normal system of linear
ordinary differential equations:
\begin{equation}
\label{eq:6} \frac{{dz_k }}{{dy}} = \sum\limits_{j = 1}^N {a_{kj}
(y)z_j} \, , \quad k = \overline{1,N} \, ,
\end{equation}
where $a_{kj}(y)$ are piecewise-continuous functions on $[y_0
,y_n]$ satisfying Lipschitz-type conditions and,
therefore, bounded, i.e. $\left| {a_{kj} (y)} \right| \le M$, with
constant $M$ and $k,j = \overline{1,N}$.

We have for system \eqref{eq:6} the next $N$ linear homogeneous
boundary conditions:
\begin{equation}
\label{eq:8}
\begin{gathered}
z_r(y_0 ) = 0 \, , \quad r = \overline{1,m} \, , \\
z_s(y_n ) = 0 \, , \quad s = \overline{m + 1,N} \, .
\end{gathered}
\end{equation}

At the discontinuity points $y = y_i$, $i = \overline{1,n - 1}$, of
functions $a_{kj}(y)$, the solution $z_k(y)$ of the system
\eqref{eq:6} should satisfy the next conjugation conditions:
\begin{equation}
\label{eq:7}
\sum\limits_{k = 1}^N d_{j,k}^{} (y_i )z_k^{(i)} (y_i )
= \sum\limits_{k = 1}^N {b_{j,k}^{} (y_i )z_k^{(i + 1)} (y_i )} \, ,
\end{equation}
$j = \overline{1,N}$, $i = \overline{1,n - 1}$, and where
$d_{j,k}(y_i)$ is the value of the function at the $i$-th interval,
and $b_{j,k}(y_i)$ is the value of the function at the $(i+1)$-th
interval.

Now, we consider a Cauchy problem for the system \eqref{eq:6} for
each interval $y_{i - 1} \le y \le y_i$, $i = \overline{1,n}$, with
the initial conditions at points $y = y_{i - 1}$. According to
Picard's Theorem, each such problem has a unique solution satisfying
the given initial conditions. Therefore, we are in conditions to
apply some well-known numerical method, for example the Runge-Kutta
method, and find the fundamental system of solutions for each $i$-th
interval, \textrm{i.e.}, to find $N$ linearly independent solutions
of the system \eqref{eq:6} in the form
\begin{equation}
\label{eq:9}
\Phi ^{(i)} (y) = \left| {\begin{array}{*{20}c}
   {\varphi _{1,1}^{(i)} (y)} & {...} & {\varphi _{1,N}^{(i)} (y)}  \\
   {...} & {...} & {...}  \\
   {\varphi _{N,1}^{(i)} (y)} & {...} & {\varphi _{N,N}^{(i)} (y)}  \\
\end{array}} \right| = \{\varphi _{k,j}^{(i)}(y) \}_{k,j} \, ,
\end{equation}
where each function $\varphi_{k,j}^{(i)}(y)$ is defined and
continuous at $[y_{i - 1} ,y_i]$, $i = \overline{1,n}$; $k$ is a
solution number and $j$ is a function number, $k,j =
\overline{1,N}$. At the points $y = y_{i - 1}$, these solutions
satisfy the conditions
\begin{equation*}
   \Phi ^{(i)} (y_{i - 1} ) = E \, , \quad i = \overline{1,n} \, ,
\end{equation*}
where $E$ is the unit matrix, or, in a more detailed way,
\begin{equation}
\label{eq:10}
\{\varphi _{k,j}^{(i-1)}(y) \}_{k,j} =
\begin{cases}
1, & k=j, \\
0, & k \neq j,
\end{cases}
\quad i = \overline{1,n}.
\end{equation}
Such fundamental system of solutions \eqref{eq:9} with initial
conditions \eqref{eq:10} is called a \emph{normal fundamental
system} \cite{c2}. We remark that for each point $y = y_{i - 1}$, $i
= \overline{1,n}$, the existence and uniqueness theorem applies. The
general solution of the system \eqref{eq:6} at the $i$-th interval
may be written with the help of the normal fundamental system of
solutions as
\begin{equation}
\label{eq:11}
z_k^{(i)} (y) = \sum\limits_{j = 1}^N
{C_j^{(i)} \varphi_{jk}^{(i)} (y)} \, ,
\end{equation}
with $C_j^{(i)}$ arbitrary constants and $y_{i - 1} \le y \le y_i$,
$k = \overline{1,N}$, $i = \overline{1,n}$. Satisfaction of the
conjugation conditions \eqref{eq:7} allows to obtain the
coefficients $C_j^{(i)}$ in the form
\begin{equation}
\label{eq:12}
C_j^{(i)}  = \sum\limits_{q = m + 1}^N {C_q^{(1)} u_{j,q}^{(i)} } \, ,
\end{equation}
with the coefficients $u_{j,q}^{(i)}$ given by recurrence as
follows:

\begin{equation}
\label{eq:coef} u_{j,q}^{(i+1)}  = \frac{\left|\begin{array}{ccccccc}
                                                 b_{1,1} & ... & b_{1,j-1} & \sum\limits_{k = 1}^N{d_{1,k}}
                                                 \sum\limits_{p = 1}^N{u_{p,q}^{(i)} \varphi_{p,k}^{(i)}(y_{i})}
                                                 & b_{1,j+1} & ... & b_{1,N} \\
                                                 ... & ... & ... & ... & ... & ... & ... \\
                                                 b_{N,1} & ... & b_{N,j-1} & \sum\limits_{k = 1}^N{d_{N,k}}
                                                 \sum\limits_{p = 1}^N{u_{p,q}^{(i)} \varphi_{p,k}^{(i)}(y_{i})}
                                                 & b_{N,j+1} & ... & b_{N,N}
                                               \end{array}
 \right|}{\left| \begin{array}{ccccc}
                                                b_{1,1} &  b_{1,2} &  ... & b_{1,N} \\
                                                ...  & ... &  ... & ... \\
                                                b_{N,1} &  b_{N,2} & ... &  b_{N,N}
                                              \end{array}
\right|}
 \, ,
\end{equation}
where  $j=\overline{1,N}$, $q=\overline{m+1,N}$ and
$i=\overline{1,n-1}$.

For application of the formulas \eqref{eq:coef} we need to know the
initial coefficients $u_{j,q}^{(1)}$ for all $j=\overline{1,N}$ and
$q=\overline{m+1,N}$. From the boundary conditions \eqref{eq:8} at
$y = y_0$ we can find
\begin{equation}
\label{eq:ufun} u_{j,q}^{(1)} =
\begin{cases}
1, & j=q, \\
0, & j \neq q,
\end{cases}
\quad \overline{j=1,N}  \, , \quad q=\overline{m+1,N} \, .
\end{equation}

Thus, the solution of the system \eqref{eq:6} for $y_{0} \leq y \leq
y_{n}$, which satisfies to the conjugation conditions \eqref{eq:7},
is determined by the formula \eqref{eq:11} subject to
\eqref{eq:12}--\eqref{eq:ufun}. According to the existence theorem,
such solution exists and it is unique at the interval
$y_{0} \leq y \leq y_{n}$.

Imposing also the initial condition \eqref{eq:8} at $y = y_n$, we
obtain a homogeneous system of $m$ linear algebraic equations in the
coefficients $C_q^{(1)}$ for all $q=\overline{m+1,N}$:
\begin{equation}
\label{eq:13} \sum\limits_{q = m + 1}^N {C_q^{(1)} } \sum\limits_{j
= 1}^N u_{j,q}^{(n)} \varphi _{j,s}^{(n)} (y_n ) = 0  \, ,
\end{equation}
where $s = \overline{m + 1,N}$.

The boundary problem \eqref{eq:6}--\eqref{eq:7} has a non-trivial
solution only at the condition of existence of non-trivial solution
of the system \eqref{eq:13}.

Thus, we conclude that a necessary and sufficient condition for the
existence of a non-trivial solution of the boundary problem
\eqref{eq:6}-\eqref{eq:8} is given by
\begin{equation}
\label{eq:14}
D = \det \left| {\sum\limits_{j = 1}^N {u_{j,q}^{(n)}
\varphi _{j,s}^{(n)} (y_n )} } \right| = 0 \, ,
\end{equation}
$s,q = \overline{m + 1,N}$. Taking into account the dependencies of
the functions $\varphi _{j,s}^{(n)}(y_n)$ on the value $\lambda$, we
determine the eigenfrequencies of the oscillation system as roots of
the equation \eqref{eq:14}.

Given a concrete system, the necessary and sufficient condition
\eqref{eq:14} provides a general method to investigate the
frequency spectrum of possible oscillations and its dependence on
different parameters of the system. In the next section we review
the application of this method for the study of different
oscillation problems.


\section{Some Applications}

Let us mention in this section some applications of the method of
normal fundamental systems of solutions for solving oscillation
problems in different technical systems.

The authors of \cite{c1} apply this method for analysis of dynamics
of some controlled machine units, which are widespread in different
modern apparatuses, machines for materials processing,
transportation, etc. The considered controlled machine units consist
of the motor, some mechanical element and the elastic connection of
significant length between them. The mathematical model of such
controlled dynamical system with distributed and discrete parameters
is a partial differential equation satisfying non-classical boundary
conditions as follows

\begin{equation}
\label{eq:paper1} GI_{p}\frac{\partial^{2}\varphi}{\partial x^{2}}-
\rho I_{p}\frac{\partial^{2}\varphi}{\partial t^{2}}+
\zeta_{1}\frac{\partial^{3}\varphi}{\partial t\partial x^{2}}=0 \, ,
\end{equation}

\begin{equation}
\label{eq:paper2} x=0 \, , \quad
J_{1}\frac{\partial^{2}\varphi}{\partial t^{2}}-
\zeta_{1}\frac{\partial^{2}\varphi}{\partial t \partial x}+ \beta
\frac{\partial \varphi}{\partial t} - GI_{p}\frac{\partial
\varphi}{\partial x}=0 \, ,
\end{equation}

\begin{equation}
\label{eq:paper3} x=L \, , \quad
J_{2}\frac{\partial^{2}\varphi}{\partial t^{2}}+
\zeta_{1}\frac{\partial^{2}\varphi}{\partial t \partial x}+ GI_{p}
\frac{\partial \varphi}{\partial x} - m_{0} \alpha_{1}\frac{\partial
\varphi}{\partial t}=0 \, ,
\end{equation}
where $\varphi(x,t)$ is the shaft arbitrary section twist; $G$ is
the shear modulus; $I_{p}$ is the polar inertia moment; $\rho$ is
the line density of the unit of the elastic link material length;
$J_{1}$ is the inertia moment of the motor rotor; $J_{2}$ is the
reduced shaft inertia moment of intermediate mechanisms
of the motor and the tool; $\zeta_{1}$ is the dissipation coefficient in the
shaft; $\beta$ is the rigidity of the motor mechanical
characteristic; $m_{0}$ is the tool mass; $\alpha_{1}$ is the
coefficient of the linear rate component of the tool resistance
moment.

The numerical method of the normal fundamental systems allows to
determine the frequencies of possible oscillations for the case of
the problem \eqref{eq:paper1}--\eqref{eq:paper3}. The influence of
different parameters on the frequencies is also analyzed by the
authors. It is found the possibility to avoid undesirable
oscillations with the help of application of different feedbacks,
which purposefully change the rigidity of the motor mechanical
characteristic in the machine unit.

In the paper \cite{c5}, we apply the numerical method of normal
fundamental systems of solutions for investigation of oscillations
in a flexible elastic system of a spacecraft of big size.

The basic mathematical model, which is used for investigations, is:

\begin{equation}
\label{eq:4f} \rho\,S  \frac{\partial^2 u}{\partial t^2} -ES
\frac{\partial^2 u}{\partial x^2} =ES\beta\, \frac{\partial^3
u}{\partial x^2 \partial t} \, ,
 \end{equation}

\begin{equation}
\label{eq:5f} x = 0 : \quad u \left( 0,t \right) =0 \, ,
 \end{equation}

\begin{multline}
\label{eq:6f} x = l : \quad mES\beta\,  \frac{\partial^2 u}{\partial
x \partial t^3} + m(d+b) \frac{\partial^3 u}{\partial t^3}
 + E S(m + \beta b) \frac{\partial^3 u}{\partial x \partial t^2}\\
+mc \frac{\partial^2 u}{\partial t^2} + E S(b + \beta\,c)
\frac{\partial^2 u}{\partial x \partial t} +ESc \frac{\partial
u}{\partial x}  =0 \, ,
 \end{multline}
where $u(x,t)$ is the displacement of a point of the bar with
abscissa $x$, $S$ the bar cross-section area, $\rho$ the unit volume
mass, $E$ the modulus of the material elasticity, $\beta$ the
coefficient of dissipation in the bar material, $b$ the damping
factor of the executive mechanism, $c$ the rigidity of the centering
spring, $d$ the feedback coefficient, $m$ the specimen mass, and $l$
is the length of the elastic bar.

The authors show that for each such dynamical system it is possible
to characterize possible oscillations using normal
fundamental functions. Moreover, it is possible to change,
purposefully, eigenvalues of the considered nonconservative boundary
problem. With this aim, different feedbacks may be applied,
changing the feedback coefficient $d$ in a desirable way. This allows
to avoid in the system an excitation of undesirable oscillations,
including self-excited ones.

These results of the considered numerical method coincide with the
ones found in \cite{c4}, obtained to the problem
\eqref{eq:4f}--\eqref{eq:6f} with a different approach: the
approximate analytical method of small parameter.

The analysis of dynamics for some aerial cable-way systems, by the
method of normal fundamental systems of solutions, may be found in
\cite{c6}. If discrete loads with masses $m_{1}$, $m_{2}$, ...,
$m_{n-1}$ move along the cable-way at a steady speed $v$, then the
mathematical model for such controlled dynamical hybrid system is
the next non-stationary boundary problem:

\begin{equation}
\label{eq:motion-disc} {\rho\frac{\partial^{2}u_{i}}{\partial
t^{2}}=T\frac{\partial^{2}u_{i}}{\partial x^{2}} \, ,}
\end{equation}

\begin{equation}
\label{bound-cond-disc} {u_{1}(0,t)=u_{n}(l,t)=0 \, ,}
\end{equation}
with conjugation conditions

\begin{equation}
\label{conjug-cond-disc1}
{u_{i+1}(l_{i},t)=u_{i}(l_{i},t)|_{x=l_{i}(t)}} \, ,
\end{equation}

\begin{equation}
\label{conjug-cond-disc2} {m_{i}(\frac{\partial^{2}u_{i}}{\partial
t^{2}}+2v\frac{\partial^{2}u_{i}}{\partial x\partial
t}}+v^{2}\frac{\partial^{2}u_{i}}{\partial
x^{2}})|_{x=l_{i}(t)}=T(\frac{\partial u_{i+1}}{\partial
x}-\frac{\partial u_{i}}{\partial x})|_{x=l_{i}(t)} \, ,
\end{equation}
where $t$ is the time variable; $x$ is the coordinate of some point
of the cable; $\rho$ is the mass per unit length; $T$ is the tension
in the cable; $l$ is the length of the cable; $m_{i}$ are the masses
of the discrete loads at arbitrary points $x=l_{i}(t)$ for all
$i=\overline{1;n-1}$; $l_{i}(t)=l_{i}+v t$, \, $l_{0}(t)=0$, \,
$t_{0} \leq t \leq t_{1}$;  $l_{i}$ is the coordinate of the
arbitrary point with the discrete mass $m_{i}$ at the initial moment
of time $t=0$, $i=\overline{1,n-1}$, $l_{0}=0$, $l_{n}=l$;
$u_{i}(x,t)$ are the functions of deviation from equilibrium
position for the point $x\in[l_{i-1},l_{i}]$ at the moment $t$,
$i=\overline{1,n}$.

For the problem \eqref{eq:motion-disc}--\eqref{conjug-cond-disc2},
which is a continuous-discrete problem, the method of the normal
fundamental functions is very efficient. It allows to find the
dependence of the dynamical characteristics of the system on the
controlled parameter, which is the speed $v$.

Finally, let us mention one open problem which may be solved by the
method of normal fundamental functions. In the paper \cite{c3}, the
approximate analytical method is used for investigations of possible
oscillations in the transmission pipelines intended for lifting
minerals from great depth.

A pipe of significant length with a big mass platform on the end is
a basic constructive element of such deep-water plants. Intensive
dynamical processes of different physical nature may be excited
under the influence of forced wave disturbances and nonlinear
hydrodynamical forces in the system. Water flow, which washes the
pipe, can continuously transport the energy into the system,
significantly exciting the system and causing dangerous oscillations.
As a result, breaking of the plant is possible.

The mathematical model of vibrations excitation, in the case when
the pipe is washed by a water flow, is a dissipative wave equation of
longitudinal vibrations:
\begin{equation}
\label{math/1} \frac{\partial^2u}{\partial
t^2}-a^2\frac{\partial^2}{\partial x^2}\left(u+\beta\frac{\partial
u}{\partial t}\right)=0,
\end{equation}
where $u(x,t)$ is the longitudinal displacement of pipe points;
$\displaystyle a^2=\frac{E}{\rho}$ ; $E$ is the modulus of
elasticity of the material; $\rho$ is the mass of material volume
unit; $\beta$ is the coefficient characterizing internal friction in
the material.

For transmission pipelines with the elastically attached pipe,
boundary conditions for the equation (\ref{math/1}) are the following:
\begin{eqnarray}
x=0,& & ES\frac{\partial}{\partial x}\left(u(0,t)+\beta\frac{\partial u(0,t)}
{\partial t}\right)=ku(0,t),\label{math/2}\\
x=L,& & M\frac{\partial^2u(L,t)}{\partial t^2}+ES\frac{\partial}
{\partial x}\left(u(L,t)+\beta\frac{\partial u(L,t)}{\partial t}\right)\nonumber\\
&&=\alpha_1\frac{\partial u(L,t)}{\partial t}-\alpha_3\left(\frac{\partial
u(L,t)}{\partial t}\right)^3,\label{math/3}
\end{eqnarray}
where $M$ is the platform mass; $S$ is the pipe cross-section area;
$k$ is the longitudinal rigidity of the elastic hanger at the point
of pipe attaching; $L$ is the pipe length; $\displaystyle
\alpha_1\frac{\partial u(L,t)}{\partial t}-\alpha_3\left(\frac{\partial
u(L,t)}{\partial t}\right)^3$ is the force of the non-linear resistance,
which appears as a result of the interaction between the body and
the environment; $\alpha_1$, $\alpha_3$ are constants characterizing
the resistant environment.

It is obvious that the method of normal fundamental functions may be
used for numerical solving the boundary problem
(\ref{math/1})--(\ref{math/3}), which was solved in the paper
\cite{c3} only approximately and with some additional assumptions.

For all the above mentioned systems, the technique for determining
complex eigenvalues of the corresponding boundary problem is
developed. With the help of the method of normal fundamental systems
of solutions, the dependencies of the oscillation frequencies on
different physical parameters of the systems are obtained. Some
possibilities of control of the frequency spectrum, in which
oscillations are possible, are found for each mentioned system.


\section*{Acknowledgments}

The authors gratefully acknowledge the support of the
\emph{Portuguese Foundation for Science and Technology} (FCT): Olena
Mul through the fellowship SFRH/BPD/14946/2004; Delfim Torres
through the R\&D unit \emph{Centre for Research in Optimization and
Control} (CEOC).


\end{document}